\documentclass{article}
\usepackage{amsfonts}
\begin{document}

\title{Logical, conditional, and classical probability}
\author{G. A. Quznetsov \\
%EndAName
gunn@chelcom.ru, lak@cgu.chel.su, gunn@mail.ru}
\date{}
\maketitle

\begin{abstract}
The propositional logic is generalized on the real numbers field. the 
logical function with all properties of the classical probability function is 
obtained. The logical analog of the Bernoulli independent tests scheme is 
constructed. The logical analog of the Large Number Law is deduced from 
properties of these functions. The logical analog of thd conditional 
probability is defined. Consistency encured by a model on a suitable variant of 
the nonstandard analysis.
\end{abstract}

\section{Introduction}
\label{intro}
There is the evident nigh affinity between the classical probability
function and the Boolean function of the classical propositional logic \cite{LYN66}. These functions are differed by the range of value, only. That is if
the range of values of the Boolean function shall be expanded from the
two-elements set $\left\{ 0;1\right\} $ to the segment $\left[ 0;1\right] $
of the real numeric axis then the logical analog of the Bernoulli Large
Number Law \cite{BER13} can be deduced from the logical axioms. These topics
is considered in this article.
\section{The classical logic}
\label{sec:1}
{\bf Definition 2.1} Sentence $\ll \Theta \gg $ is a true
sentence if and only if $\Theta $ \cite{TAR44}.

For example: sentence $\ll $it rains$\gg $ is the true sentence if and
only if it rains.

{\bf Definition 2.2} Sentence $\ll \Theta \gg $ is a false
sentence if and only if it is not that $\Theta $.

{\bf Definition 2.3} Sentences $A$ and $B$ are equal ($A=B$%
) if $A $ is true if and only if $B$ is true.

Hereinafter we use the usual notions of the classical propositional logic 
\cite{MEN63}.

{\bf Definition 2.4} Sentence $C$ is a conjunction of the
sentences $A$ and $B$ \\($C=\left( A\wedge B\right) $) if $C$ is true if and
only if $A$ is true and $B$ is true.

{\bf Definition 2.5} Sentence $C$ is a negation of the
sentence $A$ ( $C=\overline{A}$), if $C$ is true if and only if $A$ is false.

{\bf Theorem 2.1}

1) $(A\wedge A)=A$;

2) $(A\wedge B)=(B\wedge A)$;

3) $(A\wedge (B\wedge C))=((A\wedge B)\wedge C)$;

4) if $T$ is the true sentence then for every sentence $A$: $(A\wedge T)=A$;

5) if $F$ is false sentence then $\overline{F}$ is true sentence. 

{\bf Proof of the Theorem {2.1}: }From Definitions {2.1}, 
{2.2}, {2.3}, {2.4}.

{\bf Definition 2.6} Each function $\rm{g}$ with domain in 
the set of the sentences and with the range of values on the two-elements set 
$\left\{ 0;1\right\} $ is a Boolean function if:

1) $\rm{g}\left( \overline{A}\right) =1-$ ${%
\rm{g}}\left( A\right) $ for every sentence $A$;

2) $\rm{g}\left( A\wedge B\right) ={%
\rm{g}}\left( A\right) \cdot \rm{g}\left( B\right) $ for all 
sentences $A$ and $B$.

%\newpage

{\bf Definition 2.7} Set $\Im $ of the sentences is a basic set if for every element 
$A$ of this set there exist Boolean functions $\rm{g}_1$ and $\rm{g}_2$ such 
that the following conditions fulfill:

1) $\rm{g}_1\left( A\right) \neq \rm{g}_2\left( A\right) $;

2) $\rm{g}_1\left( B\right) =\rm{g}_2\left( B\right) $ for 
each element $B$ of $\Im $ such that $B\neq A$.

{\bf Definition 2.8} Set $\left[ \Im \right] $ of the sentences is 
a propositional closure of the set $\Im $ if the following
conditions fulfill:

1) if $A\in \Im $ then $A\in \left[ \Im \right ] $;

2) if $A\in \left[\Im \right]$ then $\overline{A}\in \left[ \Im \right ] $;

3) if $A\in \left[ \Im \right ]$ and $B\in \left[ \Im \right ] $ then $%
\left( A\wedge B\right) \in \left[ \Im \right ] $;

4) there do not exist other elements of $\left[ \Im \right ] $ except the 
listed by 1), 2), 3) points of this definition.

In the following text the elements of $\left[ \Im \right] $ are called as 
the $\mathit{\Im }$-sentences.

{\bf Definition 2.9} $\Im $-sentence $A$ is a tautology
if for all Boolean functions $\rm{g}$:

\[
\rm{g}(A)=1\mbox{.} 
\]

{\bf Definition 2.10} A disjunction and an
implication are defined by the usual way:

\[
\begin{array}{c}
\left( A\vee B\right) =\overline{\left( \overline{A}\wedge \overline{B}%
\right) }\mbox{,} \\ 
\left( A\Rightarrow B\right) =\overline{\left( A\wedge \overline{B}\right) }%
\mbox{.}
\end{array}
\]

By this definition and the Definitions {2.4} and {2.5}:

$\left( A\vee B\right) $ is the false sentence if and only if $A$ is the
false sentence and $B$ is the false sentence.

$\left( A\Rightarrow B\right) $ is the false sentence if and only if $A$
is the true sentence and $B$ is the false sentence.

{\bf Definition 2.11} A $\Im $-sentence is a propositional axiom 
 \cite{MEN63} if this sentence has got one some amongst the following forms:

\textbf{A1}. $\left( A\Rightarrow \left( B\Rightarrow A\right) \right) $;

\textbf{A2. }$\left( \left( A\Rightarrow \left( B\Rightarrow C\right)
\right) \Rightarrow \left( \left( A\Rightarrow B\right) \Rightarrow \left(
A\Rightarrow C\right) \right) \right) $;

\textbf{A3}. $\left( \left( \overline{B}\Rightarrow \overline{A}\right)
\Rightarrow \left( \left( \overline{B}\Rightarrow A\right) \Rightarrow
B\right) \right) $.

Let $\Im $ be some basic set. In the following text I consider $\Im $-sentences, 
only. 

{\bf Definition 2.12} Sentence $B$ is obtained from the sentences $\left( A\Rightarrow
B\right) $ and $A$ by the logic rule "modus ponens".

{\bf Definition 2.13} \cite{MEN63} Array $A_1,A_2,\ldots ,A_n$ of the
sentences is a propositional deduction of the sentence $A$ from
the hypothesis list $\Gamma $ (denote: $\Gamma \vdash A$) if $A_n=A$ and
for all numbers $l$ ($1\leq l\leq n$): $A_l$ is either the propositional
axiom or $A_l$ is obtained from some sentences $A_{l-k}$ and $A_{l-s}$ by
the modus ponens or $A_l\in \Gamma $.

{\bf Definition 2.14} A sentence is a propositional proved
sentence if this sentence is the propositional axiom or this sentence is
obtained from the propositional proved sentences by the modus ponens.

Hence, if $A$ is the propositional proved sentence then the propositional
deduction
\begin{center}
$\vdash A$
\end{center}
exists.

{\bf Theorem: 2.2} \cite{MEN63} If sentence $A$ is the propositional
proved sentence then for all Boolean function $\rm{g}$: $\rm{g}\left( A\right) =1$.

{\bf Proof of the Theorem {2.2}: }\cite{MEN63}.

{\bf Theorem: 2.3} {\bf (The completeness Theorem). }\cite{MEN63} All tautologies
are the propositional proved sentences.

{\bf Proof of the Theorem {2.3}: }\cite{MEN63}.

\section{B-functions}
\label{sec:2}

{\bf Definition 3.1} Each function $\rm{b}\left( x\right)$
with domain in the sentences set and with the range of values on the numeric 
axis segment $\left[ 0;1\right] $ is called as a B-function if 

%\begin{equation}
\[
\rm{b}\left( C\right) =1
\]
%\end{equation}

for some sentence $C$ and

\[
\rm{b}\left( A\wedge B\right) +\rm{b}\left( A\wedge 
\overline{B}\right) =\rm{b}\left( A\right)  
\]

for every sentences $A$ and $B$.

{\bf Theorem: 3.1} For each B-function $\rm{b}$:

1) for every sentences $A$ and $B$: $\rm{b}\left( A\wedge B\right)
\leq \rm{b}\left( A\right) $;

2) for every sentence $A$: if $T$ is the true sentence, then \\$\rm{b}%
\left( A\right) +\rm{b}\overline{\left( A \right )}=\rm{b}%
\left( T\right) $

3) for every sentence $A$: if $T$ is the true sentence, then $\rm{b}%
\left( A\right) \leq \rm{b}\left( T\right) $;

{\bf Proof of the Theorem {3.1}:}

1)From Definitions {3.1}.

2) From the points 4 and 2 of the Theorem {2.1}:

\[
\rm{b}\left( T\wedge A\right) +\rm{b}\left( T\wedge 
\overline{A}\right) =\rm{b}\left( A\right) +\rm{b}\left( 
\overline{A}\right) . 
\]

3) From previous point of that Theorem.

Therefore, if $T$ is the true sentence, then

\begin{equation}
\rm{b}\left( T\right) =1\mbox{.}  \label{b2}
\end{equation}

Hence, for every sentence $A$:

\begin{equation}
\rm{b}\left( A\right) +\rm{b}\left( \overline{A}\right) =1%
\mbox{.}  \label{b3}
\end{equation}

{\bf Theorem: 3.2} If sentence $D$ is the propositional proved
sentence then for all B-functions $\rm{b}$: $\rm{b}\left(
D\right) =1 $.

{\bf Proof of the Theorem {3.2}: }

If $D$ is A1 then by Definition {2.10}:

\[
\rm{b}\left( D\right) =\rm{b}\left( \overline{\left( A\wedge 
\overline{\overline{\left( B\wedge \overline{A}\right) }}\right) }\right) %
\mbox{.} 
\]

By (\ref{b3}):

\[
\rm{b}\left( D\right) =1-\rm{b}\left( A\wedge \overline{%
\overline{\left( B\wedge \overline{A}\right) }}\right) \mbox{.} 
\]

By the Definition {3.1} and the Theorem {2.1}:

\[
\begin{array}{c}
\rm{b}\left( D\right) =1-\rm{b}\left( A\right) +\rm{b}\left( A\wedge \overline{\left( B\wedge \overline{A}\right) }\right) %
\mbox{,} \\ 
\rm{b}\left( D\right) =1-\rm{b}\left( A\right) +\rm{b}\left( A\right) -\rm{b}\left( A\wedge \left( B\wedge \overline{A}%
\right) \right) \mbox{,} \\ 
\rm{b}\left( D\right) =1-\rm{b}\left( \left( A\wedge
B\right) \wedge \overline{A}\right) \mbox{,} \\ 
\rm{b}\left( D\right) =1-\rm{b}\left( A\wedge B\right) +%
\rm{b}\left( \left( A\wedge B\right) \wedge A\right) \mbox{,} \\ 
\rm{b}\left( D\right) =1-\rm{b}\left( A\wedge B\right) +%
\rm{b}\left( \left( A\wedge A\right) \wedge B\right) \mbox{,} \\ 
\rm{b}\left( D\right) =1-\rm{b}\left( A\wedge B\right) +%
\rm{b}\left( A\wedge B\right) \mbox{.}
\end{array}
\]

The proof is similar for the rest propositional axioms .

Let for all B-function $\rm{b}$: $\rm{b}(A)=1$ and $%
\rm{b}(A\Rightarrow D)=1$.

By Definition {2.10}:

\[
\rm{b}\left( A\Rightarrow D\right) =\rm{b}\left( \overline{%
A\wedge \overline{D}}\right) \mbox{.} 
\]

By (\ref{b3}):

\[
\rm{b}\left( A\Rightarrow D\right) =1-\rm{b}\left( A\wedge 
\overline{D}\right) \mbox{.} 
\]

Hence,

\[
\rm{b}\left( A\wedge \overline{D}\right) =0\mbox{.} 
\]

By Definition {3.1}:

\[
\rm{b}\left( A\wedge \overline{D}\right) =\rm{b}\left(
A\right) -\rm{b}\left( A\wedge D\right) \mbox{.} 
\]

Hence,

\[
\rm{b}\left( A\wedge D\right) =\rm{b}\left( A\right) =1%
\mbox{.} 
\]

By Definition {3.1} and the Theorem {2.1}:

\[
\rm{b}\left( A\wedge D\right) =\rm{b}\left( D\right) -{%
\rm{b}}\left( D\wedge \overline{A}\right) =1\mbox{.} 
\]

Therefore, for all B-function $\rm{b}$:

\[
\rm{b}\left( D\right) =1\mbox{.} 
\]

{\bf Theorem: 3.3}

1) If for all Boolean functions $\rm{g}$:

\[
\rm{g}\left( A\right) =1 
\]

then for all B-functions $\rm{b}$:

\[
\rm{b}\left( A\right) =1\mbox{.} 
\]

2) If for all Boolean functions $\rm{g}$:

\[
\rm{g}\left( A\right) =0 
\]

then for all B-functions $\rm{b}$:

\[
\rm{b}\left( A\right) =0\mbox{.} 
\]

{\bf Proof of the Theorem {3.3}: }

1) This just follows from the preceding Theorem and from the Theorem {2.3}.

2) If for all Boolean functions $\rm{g}$: $\rm{g}\left(
A\right) =0$, then by the Definition {2.6}: $\rm{g}\left( \overline{A}%
\right) =1$. Hence, by the point 1 of this Theorem: for all B-function ${%
\rm{b}}$: $\rm{b}\left( \overline{A}\right) =1$. By (\ref{b3}%
): $\rm{b}\left( A\right) =0$.

{\bf Theorem: 3.4} All Boolean functions are the B-functions.

Hence, the B-function is the generalization of the logic Boolean
function. Therefore, the B-function is the logic function, too.

{\bf Proof of the Theorem {3.4}: } If $C$ is {\bf A1} then 
$\rm{g}\left( C\right) =1$.

By Definition {2.6}: for all Boolean functions $\rm{g}$:

$\rm{g}\left( A\wedge B\right) +\rm{g}\left( A\wedge 
\overline{B}\right) =\rm{g}\left( A\right) \cdot \rm{g}%
\left( B\right) +{\rm{g}}\left( A\right) \cdot \left( 1-\rm{g}%
\left( B\right) \right) =\rm{g}\left( A\right) $.

{\bf Theorem: 3.5}

\[
\rm{b}\left( A\vee B\right) =\rm{b}\left( A\right) +%
\rm{b}\left( B\right) -\rm{b}\left( A\wedge B\right) \mbox{.}
\]

{\bf Definition 3.2} Sentences $A$ and $B$ are 
inconsistent sentences for the B-function $\rm{b}$ if

\[
\rm{b}\left( A\wedge B\right) =0\mbox{.} 
\]

{\bf Proof of the Theorem {3.5}: }By the Definition {2.10} and 
(\ref{b3}):

\[
\rm{b}\left( A\vee B\right) =1-\rm{b}\left( \overline{A}%
\wedge \overline{B}\right) . 
\]

By Definition {3.1}:

\[
\rm{b}\left( A\vee B\right) =1-\rm{b}\left( \overline{A}%
\right) +\rm{b}\left( \overline{A}\wedge B\right) =\rm{b}%
\left( A\right) +\rm{b}\left( B\right) -\rm{b}\left( A\wedge
B\right) \mbox{.} 
\]

{\bf Theorem: 3.6} If sentences $A$ and $B$ are the inconsistent
sentences for the B-function $\rm{b}$ then

\[
\rm{b}\left( A\vee B\right) =\rm{b}\left( A\right) +%
\rm{b}\left( B\right) \mbox{.} 
\]

{\bf Proof of the Theorem {3.6}: }This just follows from the preceding
Theorem and Definition {3.2}.

{\bf Theorem: 3.7} If $\rm{b}\left( A\wedge B\right) =\rm{b}\left( A\right)% 
\cdot \rm{b}\left( B\right) $ then $\rm{b}%
\left( A\wedge \overline{B}\right) =\rm{b}\left( A\right) \cdot {%
\rm{b}}\left( \overline{B}\right) $.

{\bf Proof of the Theorem {3.7}: }By the Definition {3.1}:

\[
\rm{b}\left( A\wedge \overline{B}\right) =\rm{b}\left(
A\right) -\rm{b}\left( A\wedge B\right) \mbox{.} 
\]

Hence,

\[
\rm{b}\left( A\wedge \overline{B}\right) =\rm{b}\left(
A\right) -\rm{b}\left( A\right) \cdot \rm{b}\left( B\right) =%
{\rm{b}}\left( A\right) \cdot \left( 1-\rm{b}\left( B\right)
\right) \mbox{.} 
\]

Hence, by (\ref{b3}):

\[
\rm{b}\left( A\wedge \overline{B}\right) =\rm{b}\left(
A\right) \cdot \rm{b}\left( \overline{B}\right) \mbox{.} 
\]

{\bf Theorem: 3.8} $\rm{b}\left( A\wedge \overline{A}\wedge
B\right) =0$.

{\bf Proof of the Theorem {3.8}: }By the Definition {3.1} and by 
the points 2 and 3 of the Theorem {2.1}:

\[
\rm{b}\left( A\wedge \overline{A}\wedge B\right) =\rm{b}%
\left( A\wedge B\right) -\rm{b}\left( A\wedge A\wedge B\right) , 
\]

hence, by the point 1 of the Theorem {2.1}:

\[
\rm{b}\left( A\wedge \overline{A}\wedge B\right) =\rm{b}%
\left( A\wedge B\right) -\rm{b}\left( A\wedge B\right) \mbox{.} 
\]

{\bf Theorem: 3.9}

\[
\mathrm{P}\left( A\wedge \left( B\vee C\right) \right) = %
\mathrm{P}\left( A\wedge B\right) +\mathrm{P}\left( A\wedge C\right) -%
\mathrm{P}\left( A\wedge B\wedge C\right) \mbox{.}
\]

{\bf Proof of the Theorem {3.9}}:

By Definition 3.1:

$\mathrm{P}\left( A\wedge \left( B\vee C\right) \right) =\mathrm{P}\left(
A\wedge \overline{\left( \overline{B}\wedge \overline{C}\right) }\right) =%
\mathrm{P}\left( A\right) -\mathrm{P}\left( A\wedge \overline{B}\wedge 
\overline{C}\right) =\mathrm{P}\left( A\right) -\mathrm{P}\left( A\wedge 
\overline{B}\right) +\mathrm{P}\left( A\wedge \overline{B}\wedge C\right) =
\mathrm{P}\left( A\wedge B\right) +\mathrm{P}\left( A\wedge C\right) -%
\mathrm{P}\left( A\wedge B\wedge C\right) $

\section{The independent tests}
\label{sec:3}

{\bf Definition 4.1} Let $st(n)$ be a function such that $st(n)$ 
has got the domain on the set of natural numbers and has got the range of 
values in the set of the $\Im $-sentences.

In this case $\Im $-sentence $A$ is a [st]-series of range $r$ with 
V- number $k$ if $A$, $r$ and $k$ fulfill to some one amongst the following 
conditions:

1) $r=1$ and $k=1$, $A=st\left( 1\right) $ or $k=0$, $A=\overline{st\left(
1\right) }$;

2) $B$ is [st]-series of range $r-1$ with V-number $k-1$ and

\[
A=\left( B\wedge st\left( r\right) \right) \mbox{,} 
\]

or $B$ is [st]-series of range $r-1$ with V-number $k$ and

\[
A=\left( B\wedge \overline{st\left( r\right) }\right) \mbox{.} 
\]

Let us denote a set of [st]-series of range $r$ with V-number $%
k$ as $[st](r,k)$.

For example, if $st\left( n\right) $ is a sentence $B_n$ then the
sentences:

$\left( B_1\wedge B_2\wedge \overline{B_3}\right) $, $\left( B_1\wedge 
\overline{B_2}\wedge B_3\right) $, $\left( \overline{B_1}\wedge B_2\wedge
B_3\right) $

are the elements of $[st](3,2)$, and \\$\left( B_1\wedge B_2\wedge \overline{%
B_3}\wedge B_4\wedge \overline{B_5}\right) \in [st](5,3)$.

{\bf Definition 4.2} Function $st(n)$ is independent for B-function $\rm{b}$ if for 
$A$: if \\$A\in $ $[st](r,r)$ then:

\[
\rm{b}\left( A\right) =\prod\limits_{n=1}^r\rm{b}\left(
st\left( n\right) \right) \mbox{.} 
\]

{\bf Definition 4.3} Let $st(n)$ be a function such that $st(n)$ 
has got the domain on the set of natural numbers and has got the range of 
values in the set of the $\Im $-sentences.

In this case sentence $A$ is [st]-disjunction of range $r$ with V-number $k$ 
(denote: $\rm{t}[st](r,k)$) if $A$ is the disjunction of all elements 
of $[st](r,k)$.

For example, if $st\left( n\right) $ is the sentence $C_n$ then:

$\left( \overline{C_1}\wedge \overline{C_2}\wedge \overline{C_3}\right) ={%
\rm{t}}[st]\left( 3,0\right) $,

$\rm{t}[st]\left( 3,1\right) =\left( \left( C_1\wedge \overline{C_2}%
\wedge \overline{C_3}\right) \vee \left( \overline{C_1}\wedge C_2\wedge 
\overline{C_3}\right) \vee \left( \overline{C_1}\wedge \overline{C_2}\wedge
C_3\right) \right) $,

$\rm{t}[st]\left( 3,2\right) =\left( \left( C_1\wedge C_2\wedge 
\overline{C_3}\right) \vee \left( \overline{C_1}\wedge C_2\wedge C_3\right)
\vee \left( C_1\wedge \overline{C_2}\wedge C_3\right) \right) $,

$\left( C_1\wedge C_2\wedge C_3\right) =\rm{t}[st]\left( 3,3\right) $. 
%\newpage

{\bf Definition 4.4} A rational number $\omega$ is called as a 
frequency of sentence $A$ in the [st]-series of $r$ independent for B-function 
$\rm{b}$ tests (designate: $\omega=\nu _r\left[ st\right] \left( A\right)$) if

1) $st(n)$ is independent for B-function $\rm{b}$,

2) for all $n$: $\rm{b}\left( st\left( n\right) \right) =\rm{b}\left( A\right) $,

3) $\rm{t}[st](r,k)$ is true and $\omega=k/r$.

{\bf Theorem: 4.1} {\bf (the J.Bernoulli formula }\cite{BER13}\textbf{)} 
If $st(n) $ is independent for B-function $\rm{b}$ and there exists 
a real number $p$ such that for all $n$: $\rm{b}\left( st\left( n\right)
\right) =p$ then

\[
\rm{b}\left( \rm{t}\left[ st\right] \left( r,k\right)
\right) =\frac{r!}{k!\cdot \left( r-k\right) !}\cdot p^k\cdot \left(
1-p\right) ^{r-k}\mbox{.} 
\]

{\bf Proof of the Theorem {4.1}: }By the Definition {4.2} and 
the Theorem {3.7}: if $B\in \left[ st\right] \left( r,k\right) $ then:

\[
\rm{b}\left( B\right) =p^k\cdot \left( 1-p\right) ^{r-k}\mbox{.} 
\]

Since $\left[ st\right] \left( r,k\right) $ contains ${r!}/\left({k!\cdot
\left( r-k\right) !}\right)$ elements then by the Theorems {3.7}, {3.8}  
and {3.6} this Theorem is fulfilled.

{\bf Definition 4.5} Let function $st(n)$ has got the domain on 
the set of the natural numbers and has got the range of values in the set of 
the $\Im $-sentences.

Let function $f(r,k,l)$ has got the domain in the set of threes of the natural 
numbers and has got the range of values in the set of the $\Im $-sentences.

In this case  $f(r,k,l)=\rm{T}[st](r,k,l)$ if

1) $f(r,k,k)=\rm{t}[st](r,k)$,

2) $f(r,k,l+1)=(f(r,k,l)\vee \rm{t}[st](r,l+1)) $.

{\bf Definition 4.6} If $a$ and $b$ are real numbers and $k-1<a\leq k$
and $l\leq b<l+1$ then $\rm{T}[st](r,a,b)=\rm{T}[st](r,k,l)$.

{\bf Theorem: 4.2} 

\[
\rm{T}[st](r,a,b)=\ll \frac ar \leq \nu _r\left[ st\right] \left(
A\right) \leq \frac br\gg \mbox{.} 
\]

{\bf Proof of the Theorem {4.2}: }By the Definition {4.6}: there 
 exist natural numbers $r $ and $k$ such that $k-1<a\leq k$ and $l\leq b<l+1$.

The recursion on $l$:

1. Let $l=k$.

In this case by the Definition {4.4}:

\[
\rm{T}[st](r,k,k)=\rm{t}[st](r,k)=\ll \nu _r\left[ st\right]
\left( A\right) =\frac kr\gg \mbox{.} 
\]

2. Let $n$ be any natural number.

The recursive assumption: Let

\[
\rm{T}[st](r,k,k+n)=\ll \frac kr\leq \nu _r\left[ st\right] \left(
A\right) \leq \frac {k+n}r\gg \mbox{.} 
\]

By the Definition {4.5}:

\[
\rm{T}[st](r,k,k+n+1)=(\rm{T}[st](r,k,k+n)\vee \rm{t}%
[st](r,k+n+1))\mbox{.} 
\]

By the recursive assumption and by the Definition {4.4}:

\[
\rm{T}[st](r,k,k+n+1)= 
\]

\[
=(\ll \frac kr\leq \nu _r\left[ st\right] \left( A\right) \leq \frac
{k+n}r\gg \vee \ll \nu _r\left[ st\right] \left( A\right) =\frac {k+n+1}r\gg)%
\mbox{.} 
\]

Hence, by the Definition {2.10}:

\[
\rm{T}[st](r,k,k+n+1)=\ll \frac kr\leq \nu _r\left[ st\right] \left(
A\right) \leq \frac {k+n+1}r\gg \mbox{.} 
\]

{\bf Theorem: 4.3} If $st(n)$ is independent for B-function ${%
\rm{b}}$ and there exists a real number $p$ such that ${%
\rm{b}}\left( st\left( n\right) \right) =p$ for all $n$ then

\[
\rm{b}\left( \rm{T}[st](r,a,b)\right) =\sum_{a\leq k\leq
b}\frac {r!}{k!\cdot \left( r-k\right) !}\cdot p^k\cdot \left( 1-p\right)
^{r-k}\mbox{.} 
\]

{\bf Proof of the Theorem {4.3}: }This is the consequence from the 
Theorem {4.1} by the Theorem {3.6}.

{\bf Theorem: 4.4} If $st(n)$ is independent for the B-function ${%
\rm{b}}$ and there exists a real number $p$ such that ${%
\rm{b}}\left( st\left( n\right) \right) =p$ for all $n$ then

\[
\rm{b}\left( \rm{T}[st](r,r\cdot \left( p-\varepsilon
\right) ,r\cdot \left( p+\varepsilon \right) )\right) \geq 1-\frac{p\cdot
\left( 1-p\right) }{r\cdot \varepsilon ^2} 
\]

for every positive real number $\varepsilon $.

{\bf Proof of the Theorem {4.4}:} Because

\[
\sum_{k=0}^r\left( k-r\cdot p\right) ^2\cdot \frac{r!}{k!\cdot \left(
r-k\right) !}\cdot p^k\cdot \left( 1-p\right) ^{r-k}=r\cdot p\cdot \left(
1-p\right) 
\]

then if

\[
J=\left\{ k\in \mathbf{N}|0\leq k\leq r\cdot \left( p-\varepsilon \right)
\right\} \cap \left\{ k\in \mathbf{N}|r\cdot \left( p+\varepsilon \right)
\leq k\leq r\right\} 
\]

then

\[
\sum_{k\in J}\frac {r!}{k!\cdot \left( r-k\right) !}\cdot p^k\cdot \left(
1-p\right) ^{r-k}\leq \frac {p\cdot \left( 1-p\right) }{r\cdot \varepsilon
^2}\mbox{.} 
\]

Hence, by (\ref{b3}) this Theorem is fulfilled.

Hence

\begin{equation}
\lim\limits_{r\rightarrow \infty }\rm{b}\left( \rm{T}[st](r,r%
\cdot \left( p-\varepsilon \right) ,r\cdot \left( p+\varepsilon \right) )%
\right) =1  \label{uh}
\end{equation}

for all tiny positive numbers $\varepsilon $.

\section{The logic probability function}
\label{sec:4}

{\bf Definition 5.1} B-function $\mathrm{P}$ is $P$-function
if for every $\Im $-sentence $\ll \Theta \gg$: \\If $\mathrm{P \left( \ll
\Theta \gg \right) = 1}$ then $\ll \Theta \gg$ is true sentence.

Hence from Theorem {4.2} and (\ref{uh}): if $\rm{b}$ is a $P$-function
then the sentence

\[
\ll \left( p-\varepsilon \right) \leq \nu _r\left[ st\right] \left( A\right)
\leq \left( p+\varepsilon \right) \gg 
\]

is almost true sentence for large $r$ and for all tiny $\varepsilon $.
Therefore, it is almost truely that

\[
\nu _r\left[ st\right] \left( A\right) =p 
\]

for large $r$.

Therefore, it is almost true that

\[
\rm{b}\left( A\right) =\nu _r\left[ st\right] \left( A\right) 
\]

for large $r$.

Therefore, the function, defined by the Definition {5.1} has got the
statistical meaning. That is why I'm call such function as the
logic probability function.

\section{Conditional probability}

{\bf Definition 6.1:} {\it Conditional probability} $B$ for $C$ is the 
following function: 

\begin{equation}
\mathfrak{b}\left( B/C\right) \stackrel{def}{=}\frac{\mathfrak{b}\left( C\wedge B\right) 
}{\mathfrak{b}\left( C\right) }\mbox{.}\label{CP}
\end{equation}

{\bf Theorem 6.1} The conditional probability function is a B-function.

{\bf Proof of Theorem 6.1} From Definition 6.1:

\begin{center}
$\mathfrak{b}\left( C/C\right) =\frac{\mathfrak{b}\left( C\wedge C\right) }{\mathfrak{b}%
\left( C\right) }$.
\end{center}

Hence by point 1 of Theorem 2.1:

\begin{center}
$\mathfrak{b}\left( C/C\right) =\frac{\mathfrak{b}\left( C\right) }{\mathfrak{b}\left(
C\right) }=1$.
\end{center}

Form Definition 6.1:

\begin{center}
$\mathfrak{b}\left( \left( A\wedge B\right) /C\right) +\mathfrak{b}\left( \left(
A\wedge \left( \neg B\right) \right) /C\right) =\frac{\mathfrak{b}\left( C\wedge \left(
A\wedge B\right) \right) }{\mathfrak{b}\left( C\right) }+\frac{\mathfrak{b}\left(
C\wedge \left( A\wedge \left( \neg B\right) \right) \right) }{\mathfrak{b}\left( C\right) }
$.
\end{center}

Hence:

\begin{center}
$\mathfrak{b}\left( \left( A\wedge B\right) /C\right) +\mathfrak{b}\left( \left(
A\wedge \left( \neg B\right) \right) /C\right) =\frac{\mathfrak{b}\left( C\wedge \left(
A\wedge B\right) \right) +\mathfrak{b}\left( C\wedge \left( A\wedge \left( \neg B\right)
\right) \right) }{\mathfrak{b}\left( C\right) }$.
\end{center}

By point 3 of Theorem 2.1:

\begin{center}
$\mathfrak{b}\left( \left( A\wedge B\right) /C\right) +\mathfrak{b}\left( \left(
A\wedge \left( \neg B\right) \right) /C\right) =\frac{\mathfrak{b}\left( \left(
C\wedge A\right) \wedge B\right) +\mathfrak{b}\left( \left( C\wedge A\right) \wedge \left( \neg
B\right) \right) }{\mathfrak{b}\left( C\right) }$.
\end{center}

Hence by Definition 3.1:

\begin{center}
$\mathfrak{b}\left( \left( A\wedge B\right) /C\right) +\mathfrak{b}\left( \left(
A\wedge \left( \neg B\right) \right) /C\right) =\frac{\mathfrak{b}\left( C\wedge A\right) 
}{\mathfrak{b}\left( C\right) }$.
\end{center}

Hence by Definition 6.1:

\begin{center}
$\mathfrak{b}\left( \left( A\wedge B\right) /C\right) +\mathfrak{b}\left( \left(
A\wedge \left( \neg B\right) \right) /C\right) =\mathfrak{b}\left( A/C\right) $ 
$_{\bf \Box }$
\end{center}

\section{Classical probability}

Let $\mathrm{P}$ be $P$-function.

{\bf Definition 7.1} $\left\{ B_1,B_2,\ldots ,B_n\right\} $ is called as {\it 
complete set} if the following conditions are fulfilled:

1. if $k\neq s$ then $\left( B_k\wedge  B_s\right)$ is a false sentence;

2. $\left( B_1\vee B_2\vee \ldots \vee B_n\right)$ is a true sentence.

{\bf Definition 7.2} $B$ is favorable for $A$ if $\left(B\wedge\overline{A}%
\right)$ is a false sentence, and $B$ is unfavorable for $A$ if $\left(B\wedge% 
A\right)$ is a false sentence.

Let 

1. $\left\{ B_1,B_2,\ldots ,B_n\right\}$ be complete set;

2. for $k\in \left\{ 1,2,\ldots ,n\right\} $ and $s\in \left\{ 1,2,\ldots ,n%
\right\} $: $\mathrm{P}\left( B_k\right) =\mathrm{P}\left( B_s\right) $;

3. if $1\leq k\leq m$ then $B_k$ is favorable for $A$, and if $m+1\leq s\leq n$
then  $B_s$ is unfavorable for $A$.

In that case from point 5 of Theorem 2.1 and from (\ref{b2}) and (\ref{b3}):

\[
\mathrm{P}\left( \overline{A}\wedge B_k\right) = 0
\]

for $k\in \left\{ 1,2,\ldots ,m\right\} $ and 

\[
\mathrm{P}\left( A\wedge B_s\right) =0 
\]

for  $s\in \left\{ m+1,m+2,\ldots ,n\right\} $.

Hence from Definition 3.1:

\[
\mathrm{P}\left( A\wedge B_k\right) =\mathrm{P}\left( B_k\right)
\]

for $k\in \left\{ 1,2,\ldots ,n\right\} $.

By point 4 of Theorem 2.1:

\[
A=\left( A\wedge \left( B_1\vee B_2\vee \ldots \vee B_m\vee B_{m+1}\ldots
\vee B_n\right) \right) \mbox{.}
\]

Hence by Theorem 3.9:

$\mathrm{P}\left( A\right) =\mathrm{P}\left( A\wedge B_1\right) +\mathrm{P}%
\left( A\wedge B_2\right) +\ldots +$

$+\mathrm{P}\left( A\wedge B_m\right) +\mathrm{P}\left( A\wedge
B_{m+1}\right) +\ldots +\mathrm{P}\left( A\wedge B_n\right) =$

$=\mathrm{P}\left( B_1\right) +\mathrm{P}\left( B_2\right) +\ldots +\mathrm{P%
}\left( B_m\right) $.

Therefore

\[
\mathrm{P}\left( A\right) =\frac mn\mbox{.}
\]

\section{Conclusion}
\label{sect:concl}

The logic probability function is the extension of the logic B-function.
Therefore, \textbf{the probability is some generalization of the classic
propositional logic.} That is the probability is the logic of events such 
that these events do not happen, yet.

\section{Appendix. Consistency}

\subsection{THE NONSTANDARD NUMBERS} 

Let us consider the set ${\bf N}$ of natural numbers. 

{\bf Definition A.1:} The $n${\it -part-set} ${\bf S}$ of ${\bf N}$ is 
defined recursively as follows: 

1) ${\bf S}_1=\left\{ 1\right\} $; 

2) ${\bf S}_{\left( n+1\right) }={\bf S}_n\cup \left\{ n+1\right\} $. 

{\bf Definition A.2: }If ${\bf S}_n$ is the $n$-part-set of ${\bf N}$ and $% 
{\bf A}\subseteq {\bf N}$ then $\left\| {\bf A}\cap {\bf S}_n\right\| $ is 
the quantity elements of the set ${\bf A}\cap {\bf S}_n$, and if 

\[ 
\varpi _n\left( {\bf A}\right) =\frac{\left\| {\bf A}\cap {\bf S}_n\right\| }% 
n\mbox{,} 
\] 

then $\varpi _n\left( {\bf A}\right) $ is {\it the frequency} of the set $% 
{\bf A}$ on the $n$-part-set ${\bf S}_n$. 

{\bf Theorem A.1:} 

1) $\varpi _n({\bf N})=1$; 

2) $\varpi _n(\emptyset )=0$; 

3) $\varpi _n({\bf A})+\varpi _n({\bf N}-{\bf A})=1$; 

4) $\varpi _n({\bf A}\cap {\bf B})+\varpi _n({\bf A}\cap ({\bf N}-{\bf B}% 
))=\varpi _n({\bf A})$. 

{\bf Proof of the Theorem A.1:} From Definitions A.1 and A.2. 

{\bf Definition A.3: }If ''$\lim $'' is the Cauchy-Weierstrass ''limit'' 
then let us denote: 

\[ 
{\bf \Phi ix=}\left\{ {\bf A}\subseteq {\bf N}|\lim_{n\rightarrow \infty 
}\varpi _n({\bf A})=1\right\} \mbox{.} 
\] 

{\bf Theorem A.2: }${\bf \Phi ix}$ is the filter \cite{DVS}, i.e.: 

1) ${\bf N}\in {\bf \Phi ix}$, 

2) $\emptyset \notin {\bf \Phi ix}$, 

3) if ${\bf A}\in {\bf \Phi ix}$ and ${\bf B}\in {\bf \Phi ix}$ then $({\bf A% 
}\cap {\bf B})\in {\bf \Phi ix}$ ; 

4) if ${\bf A}\in {\bf \Phi ix}$ and ${\bf A}\subseteq {\bf B}$ then ${\bf B}% 
\in {\bf \Phi ix}$. 

{\bf Proof of the Theorem A.2:} From the point 3 of Theorem A.1: 

\[ 
\lim_{n\rightarrow \infty }\varpi _n({\bf N}-{\bf B})=0\mbox{.} 
\] 

From the point 4 of Theorem A.1: 

\[ 
\varpi _n({\bf A}\cap ({\bf N}-{\bf B}))\leq \varpi _n({\bf N}-{\bf B})% 
\mbox{.} 
\] 

Hence, 

\[ 
\lim_{n\rightarrow \infty }\varpi _n\left( {\bf A}\cap ({\bf N}-{\bf B}% 
)\right) =0\mbox{.} 
\] 

Hence, 

\[ 
\lim_{n\rightarrow \infty }\varpi _n\left( {\bf A}\cap {\bf B}\right) 
=\lim_{n\rightarrow \infty }\varpi _n({\bf A})\mbox{.} 
\] 

In the following text we shall adopt to our topics the definitions and the 
proofs of the Robinson Nonstandard Analysis \cite{DVS2}: 

{\bf Definition A.4:} The sequences of the real numbers $\left\langle 
r_n\right\rangle $ and $\left\langle s_n\right\rangle $ are {\it Q-equivalent% 
} (denote: $\left\langle r_n\right\rangle \sim \left\langle s_n\right\rangle 
$) if 

\[ 
\left\{ n\in {\bf N}|r_n=s_n\right\} \in {\bf \Phi ix}\mbox{.} 
\] 

{\bf Theorem A.3:} If ${\bf r}$,${\bf s}$,${\bf u}$ are the sequences of the 
real numbers then 

1) ${\bf r}\sim {\bf r}$, 

2) if ${\bf r}\sim {\bf s}$ then ${\bf s}\sim {\bf r}$; 

3) if ${\bf r}\sim {\bf s}$ and ${\bf s}\sim {\bf u}$ then ${\bf r}\sim {\bf % 
u}$.

{\bf Proof of the Theorem A.3:} By Definition A.4 from the Theorem A.2 is 
obvious. 
 
{\bf Definition A.5:} {\it The Q-number} is the set of the Q-equivalent 
sequences of the real numbers, i.e. if $\widetilde{a}$ is the Q-number and $% 
{\bf r}\in \widetilde{a}$ and ${\bf s}\in \widetilde{a}$, then ${\bf r}\sim 
{\bf s};$ and if ${\bf r}\in \widetilde{a}$ and ${\bf r}\sim {\bf s}$ then $% 
{\bf s}\in \widetilde{a}$. 

{\bf Definition A.6:} The Q-number $\widetilde{a}$ is {\it the standard 
Q-number} $a$ if $a$ is some real number and the sequence $\left\langle 
r_n\right\rangle $ exists, for which: $\left\langle r_n\right\rangle \in 
\widetilde{a}$ and 

\[ 
\left\{ n\in {\bf N}|r_n=a\right\} \in {\bf \Phi ix}\mbox{.} 
\] 

{\bf Definition A.7:} The Q-numbers $\widetilde{a}$ and $\widetilde{b}$ are 
{\it the equal Q-numbers} (denote: $\widetilde{a}=\widetilde{b}$) if a $% 
\widetilde{a}\subseteq \widetilde{b}$ and $\widetilde{b}\subseteq \widetilde{% 
a}$. 

{\bf Theorem A.4: }Let $\mathfrak{f}(x,y,z)$ be a function, which has got the 
domain in ${\bf R}\times {\bf R}\times {\bf R}$, has got the range of values 
in ${\bf R}$ (${\bf R}$ is the real numbers set). 

Let $\left\langle y_{1,n}\right\rangle $ , $\left\langle 
y_{2,n}\right\rangle $ , $\left\langle y_{3,n}\right\rangle $ , $% 
\left\langle z_{1,n}\right\rangle $ , $\left\langle z_{2,n}\right\rangle $ , 
$\left\langle z_{3,n}\right\rangle $ be any sequences of real numbers. 

In this case if $\left\langle z_{i,n}\right\rangle \sim \left\langle 
y_{i,n}\right\rangle $ then $\left\langle \mathfrak{f}(y_{1,n},y_{2,n},y_{3,n})% 
\right\rangle \sim \left\langle \mathfrak{f}(z_{1,n},z_{2,n},z_{3,n})\right% 
\rangle $.

{\bf Proof of the Theorem A.4:} Let us denote: 

if $k=1$ or $k=2$ or $k=3$ then 

\[ 
{\bf A}_k=\left\{ n\in {\bf N}|y_{k,n}=z_{k,n}\right\} \mbox{.} 
\] 

In this case by Definition A.4 for all $k$: 

\[ 
{\bf A}_k\in {\bf \Phi ix}\mbox{.} 
\] 

Because 

\[ 
\left( {\bf A}_1\cap {\bf A}_2\cap {\bf A}_3\right) \subseteq \left\{ n\in 
{\bf N}|{\mathfrak f}(y_{1,n},y_{2,n},y_{3,n})={\mathfrak f}% 
(z_{1,n},z_{2,n},z_{3,n})\right\} \mbox{,} 
\] 

then by Theorem A.2: 

\[ 
\left\{ n\in {\bf N}|{\mathfrak f}(y_{1,n},y_{2,n},y_{3,n})={\mathfrak f}% 
(z_{1,n},z_{2,n},z_{3,n})\right\} \in {\bf \Phi ix}\mbox{.} 
\] 
 
{\bf Definition A.8:} Let us denote: $Q{\bf R}$ is the set of the Q-numbers. 
~ 

{\bf Definition A.9: }The function $\widetilde{\mathfrak{f}}$, which has got 
the domain in $Q{\bf R}\times Q{\bf R}\times Q{\bf R}$, has got the range of 
values in $Q{\bf R}$, is {\it the Q-extension of the function} $\mathfrak{f}$, 
which has got the domain in ${\bf R}\times {\bf R}\times {\bf R}$, has got 
the range of values in ${\bf R}$, if the following condition is accomplished: 

Let $\left\langle x_n\right\rangle $ ,$\left\langle y_n\right\rangle $ ,$% 
\left\langle z_n\right\rangle $ be any sequences of real numbers. In this 
case: if 

$\left\langle x_n\right\rangle \in \widetilde{x}$, $\left\langle 
y_n\right\rangle \in \widetilde{y}$, $\left\langle z_n\right\rangle \in 
\widetilde{z}$, $\widetilde{u}=\widetilde{\mathfrak{f}}\left( \widetilde{x},% 
\widetilde{y},\widetilde{z}\right) $, 

then 

$\left\langle \mathfrak{f}\left( x_n,y_n,z_n\right) \right\rangle \in 
\widetilde{u}$. 

{\bf Theorem A.5:} For all functions $\mathfrak{f}$, which have the domain in $% 
{\bf R}\times {\bf R}\times {\bf R}$, have the range of values in ${\bf R}$, 
and for all real numbers $a$, $b$, $c$, $d$: if $\widetilde{\mathfrak{f}}$ is 
the Q-extension of $\mathfrak{f}$; $\widetilde{a}$, $\widetilde{b}$, $% 
\widetilde{c}$, $\widetilde{d}$ are standard Q-numbers $a$, $b$, $c$, $d$, 
then: 

if $d=\mathfrak{f}(a,b,c)$ then $\widetilde{d}=\widetilde{\mathfrak{f}}(\widetilde{% 
a},\widetilde{b},\widetilde{c})$ and vice versa. 

{\bf Proof of the Theorem A.5:} If $\left\langle r_n\right\rangle \in 
\widetilde{a}$, $\left\langle s_n\right\rangle \in \widetilde{b}$, $% 
\left\langle u_n\right\rangle \in \widetilde{c}$, $\left\langle {\mathfrak t}% 
_n\right\rangle \in \widetilde{d}$ then by Definition A.6: 

\[ 
\begin{array}{c} 
\left\{ n\in {\bf N}|r_n=a\right\} \in {\bf \Phi ix}\mbox{,} \\ 
\left\{ n\in {\bf N}|s_n=b\right\} \in {\bf \Phi ix}\mbox{,} \\ 
\left\{ n\in {\bf N}|u_n=c\right\} \in {\bf \Phi ix}\mbox{,} \\ 
\left\{ n\in {\bf N}|t_n=d\right\} \in {\bf \Phi ix}\mbox{.} 
\end{array} 
\] 

1) Let $d={\mathfrak f}(a,b,c)$. 

In this case by Theorem A.2: 

\[ 
\left\{ n\in {\bf N}|t_n={\mathfrak f}(r_n,s_n,u_n)\right\} \in {\bf \Phi ix}% 
\mbox{.} 
\] 

Hence, by Definition A.4: 

\[ 
\left\langle t_n\right\rangle \sim \left\langle {\mathfrak f}(r_n,s_n,u_n)\right% 
\rangle \mbox{.} 
\] 

Therefore by Definition A.5: 

\[ 
\left\langle {\mathfrak f}(r_n,s_n,u_n)\right\rangle \in \widetilde{d}\mbox{.} 
\] 

Hence, by Definition A.9: 

\[ 
\widetilde{d}=\widetilde{{\mathfrak f}}(\widetilde{a},\widetilde{b},\widetilde{c}% 
)\mbox{.} 
\] 

2) Let $\widetilde{d}=\widetilde{{\mathfrak f}}(\widetilde{a},\widetilde{b},% 
\widetilde{c})$. 

In this case by Definition A.9: 

\[ 
\left\langle {\mathfrak f}(r_n,s_n,u_n)\right\rangle \in \widetilde{d}\mbox{.} 
\] 

Hence, by Definition A.5: 

\[ 
\left\langle t_n\right\rangle \sim \left\langle {\mathfrak f}(r_n,s_n,u_n)\right% 
\rangle \mbox{.} 
\] 

Therefore, by Definition A.4: 

\[ 
\left\{ n\in {\bf N}|t_n={\mathfrak f}(r_n,s_n,u_n)\right\} \in {\bf \Phi ix}% 
\mbox{.} 
\] 

Hence, by the Theorem A.2: 

\[ 
\left\{ n\in {\bf N}|t_n={\mathfrak f}(r_n,s_n,u_n),r_n=a,s_n=b,u_n=c,t_n=d% 
\right\} \in {\bf \Phi ix}\mbox{.} 
\] 

Hence, since this set does not empty, then 

\[ 
d={\mathfrak f}(a,b,c)\mbox{.} 
\] 
 
By this Theorem: if $\widetilde{\mathfrak{f}}$ is the Q-extension of the 
function $\mathfrak{f}$ then the expression ''$\widetilde{\mathfrak{f}}(\widetilde{% 
x},\widetilde{y},\widetilde{z})$'' will be denoted as ''$\mathfrak{f}(% 
\widetilde{x},\widetilde{y},\widetilde{z})$'' and if $\widetilde{u}$ is the 
standard Q-number then the expression ''$\widetilde{u}$'' will be denoted as 
''$u$''. 

{\bf Theorem A.6:} If for all real numbers $a$, $b$, $c$: 

\[ 
\varphi (a,b,c)=\psi (a,b,c) 
\] 

then for all Q-numbers $\widetilde{x}$, $\widetilde{y}$, $\widetilde{z}$: 

\[ 
\varphi (\widetilde{x},\widetilde{y},\widetilde{z})=\psi (\widetilde{x},% 
\widetilde{y},\widetilde{z})\mbox{.} 
\]

{\bf Proof of the Theorem A.6:} If $\left\langle x_n\right\rangle \in 
\widetilde{x}$, $\left\langle y_n\right\rangle \in \widetilde{y}$, $% 
\left\langle z_n\right\rangle \in \widetilde{z}$, $\widetilde{u}=\varphi (% 
\widetilde{x},\widetilde{y},\widetilde{z})$, then by Definition A.9: $% 
\left\langle \varphi (x_n,y_n,z_n)\right\rangle \in \widetilde{u}$. 

Because $\varphi (x_n,y_n,z_n)=\psi (x_n,y_n,z_n)$ then $\left\langle \psi 
(x_n,y_n,z_n)\right\rangle \in \widetilde{u}$. 

If $\widetilde{v}=\psi (\widetilde{x},\widetilde{y},\widetilde{z})$ then by 
Definition A.9: $\left\langle \psi (x_n,y_n,z_n)\right\rangle \in \widetilde{% 
v}$, too. 

Therefore, for all sequences $\left\langle t_n\right\rangle $ of real 
numbers: if $\left\langle t_n\right\rangle \in \widetilde{u}$ then by 
Definition A.5: $\left\langle t_n\right\rangle \sim \left\langle \psi 
(x_n,y_n,z_n)\right\rangle $. 

Hence, $\left\langle t_n\right\rangle \in \widetilde{v}$; and if $% 
\left\langle t_n\right\rangle \in \widetilde{v}$ then $\left\langle 
t_n\right\rangle \sim \left\langle \varphi (x_n,y_n,z_n)\right\rangle $; 
hence, $\left\langle t_n\right\rangle \in \widetilde{u}$. 

Therefore, $\widetilde{u}=\widetilde{v}$. 
 
{\bf Theorem A.7:} If for all real numbers $a$, $b$, $c$: 

\[ 
\mathfrak{f}\left( a,\varphi (b,c)\right) =\psi (a,b,c) 
\] 

then for all Q-numbers $\widetilde{x}$, $\widetilde{y}$, $\widetilde{z}$: 

\[ 
\mathfrak{f}\left( \widetilde{x},\varphi (\widetilde{y},\widetilde{z})\right) 
=\psi (\widetilde{x},\widetilde{y},\widetilde{z})\mbox{.} 
\] 

{\bf Consequences from Theorems A.6 and A.7:} \cite{DVS3}: For all Q-numbers 
$\widetilde{x}$, $\widetilde{y}$, $\widetilde{z}$: 

${\bf \Phi }${\bf 1:} $(\widetilde{x}+\widetilde{y})=(\widetilde{y}+% 
\widetilde{x})$, 

${\bf \Phi }${\bf 2:} $(\widetilde{x}+(\widetilde{y}+\widetilde{z}))=((% 
\widetilde{x}+\widetilde{y})+\widetilde{z})$, 

${\bf \Phi }${\bf 3:} $(\widetilde{x}+0)=\widetilde{x}$, 

${\bf \Phi }${\bf 5:} $(\widetilde{x}\cdot \widetilde{y})=(\widetilde{y}% 
\cdot \widetilde{x})$, 

${\bf \Phi }${\bf 6:} $(\widetilde{x}\cdot (\widetilde{y}\cdot \widetilde{z}% 
))=((\widetilde{x}\cdot \widetilde{y})\cdot \widetilde{z})$, 

${\bf \Phi 7}${\bf : }$(\widetilde{x}\cdot 1)=\widetilde{x}$, 

${\bf \Phi }${\bf 10:} $(\widetilde{x}\cdot (\widetilde{y}+\widetilde{z}))=((% 
\widetilde{x}\cdot \widetilde{y})+(\widetilde{x}\cdot \widetilde{z}))$. 

{\bf Proof of the Theorem A.7:} Let $\left\langle w_n\right\rangle \in 
\widetilde{w}$, ${\mathfrak f}(\widetilde{x},\widetilde{w})=\widetilde{u}$, $% 
\left\langle x_n\right\rangle \in \widetilde{x}$, $\left\langle 
y_n\right\rangle \in \widetilde{y}$, $\left\langle z_n\right\rangle \in 
\widetilde{z}$, $\varphi (\widetilde{y},\widetilde{z})=\widetilde{w}$, $\psi 
(\widetilde{x},\widetilde{y},\widetilde{z})=\widetilde{v}$. 

By the condition of this Theorem: ${\mathfrak f}(x_n,\varphi (y_n,z_n))=\psi 
(x_n,y_n,z_n)$. 

By Definition A.9: $\left\langle \psi (x_n,y_n,z_n)\right\rangle \in 
\widetilde{v}$, $\left\langle \varphi (x_n,y_n)\right\rangle \in \widetilde{w% 
}$, $\left\langle {\mathfrak f}(x_n,w_n)\right\rangle \in \widetilde{u}$. 

For all sequences $\left\langle t_n\right\rangle $ of real numbers: 

1) If $\left\langle t_n\right\rangle \in \widetilde{v}$ then by Definition 
A.5: $\left\langle t_n\right\rangle \sim \left\langle \psi 
(x_n,y_n,z_n)\right\rangle $. 

Hence $\left\langle t_n\right\rangle \sim \left\langle {\mathfrak f}(x_n,\varphi 
(y_n,z_n))\right\rangle $. 

Therefore, by Definition A.4: 

\[ 
\left\{ n\in {\bf N}|t_n={\mathfrak f}(x_n,\varphi \left( y_n,z_n\right) 
)\right\} \in {\bf \Phi ix} 
\] 

and 

\[ 
\left\{ n\in {\bf N}|w_n=\varphi \left( y_n,z_n\right) \right\} \in {\bf % 
\Phi ix}\mbox{.} 
\] 

Hence, by Theorem A.2: 

\[ 
\left\{ n\in {\bf N}|t_n={\mathfrak f}(x_n,w_n)\right\} \in {\bf \Phi ix}\mbox{.} 
\] 

Hence, by Definition A.4: 

\[ 
\left\langle t_n\right\rangle \sim \left\langle {\mathfrak f}(x_n,w_n)\right% 
\rangle \mbox{.} 
\] 

Therefore, by Definition A.5: $\left\langle t_n\right\rangle \in \widetilde{u% 
}$. 

2) If $\left\langle t_n\right\rangle \in \widetilde{u}$ then by Definition 
A.5: $\left\langle t_n\right\rangle \sim \left\langle {\mathfrak f}% 
(x_n,w_n)\right\rangle $. 

Because $\left\langle w_n\right\rangle \sim \left\langle \varphi 
(y_n,z_n)\right\rangle $ then by Definition A.4: 

\[ 
\left\{ n\in {\bf N}|t_n={\mathfrak f}(x_n,w_n)\right\} \in {\bf \Phi ix}\mbox{,} 
\] 

\[ 
\left\{ n\in {\bf N}|w_n=\varphi \left( y_n,z_n\right) \right\} \in {\bf % 
\Phi ix}\mbox{.} 
\] 

Therefore, by Theorem A.2: 

\[ 
\left\{ n\in {\bf N}|t_n={\mathfrak f}(x_n,\varphi \left( y_n,z_n\right) 
)\right\} \in {\bf \Phi ix}\mbox{.} 
\] 

Hence, by Definition A.4: 

\[ 
\left\langle t_n\right\rangle \sim \left\langle {\mathfrak f}(x_n,\varphi 
(y_n,z_n))\right\rangle \mbox{.} 
\] 

Therefore, 

\[ 
\left\langle t_n\right\rangle \sim \left\langle \psi 
(x_n,y_n,z_n)\right\rangle \mbox{.} 
\] 

Hence, by Definition A.5: $\left\langle t_n\right\rangle \in \widetilde{v}$. 

From above and from 1) by Definition A.7: $\widetilde{u}=\widetilde{v}$. 

{\bf Theorem A.8: }${\bf \Phi }${\bf 4:} For every Q-number $\widetilde{x}$ 
the Q-number $\widetilde{y}$ exists, for which: 

$(\widetilde{x}+\widetilde{y})=0$. 

{\bf Proof of the Theorem A.8: }If $\left\langle x_n\right\rangle \in 
\widetilde{x}$ then $\widetilde{y}$ is the Q-number, which contains $% 
\left\langle -x_n\right\rangle $. 

{\bf Theorem A.9: }${\bf \Phi 9}${\bf :} There is not that $0=1$. 

{\bf Proof of the Theorem A.9:} is obvious from Definition A.6 and 
Definition A.7. 

{\bf Definition A.10:} The Q-number $\widetilde{x}$ is {\it Q-less} than the 
Q-number $\widetilde{y}$ (denote: $\widetilde{x}<\widetilde{y}$) if the 
sequences $\left\langle x_n\right\rangle $ and $\left\langle 
y_n\right\rangle $ of real numbers exist, for which: $\left\langle 
x_n\right\rangle \in \widetilde{x}$, $\left\langle y_n\right\rangle \in 
\widetilde{y}$ and 

\[ 
\left\{ n\in {\bf N}|x_n<y_n\right\} \in {\bf \Phi ix}\mbox{.} 
\] 

{\bf Theorem A.10:} For all Q-numbers $\widetilde{x}$, $\widetilde{y}$, $% 
\widetilde{z}$: \cite{DVS4} 

${\bf \Omega 1}$: there is not that $\widetilde{x}<\widetilde{x}$; 

${\bf \Omega 2}$: if $\widetilde{x}<\widetilde{y}$ and $\widetilde{y}<% 
\widetilde{z}$ then $\widetilde{x}<\widetilde{z}$; 

${\bf \Omega 4}$: if $\widetilde{x}<\widetilde{y}$ then $(\widetilde{x}+% 
\widetilde{z})<(\widetilde{y}+\widetilde{z})$; 

${\bf \Omega 5}$: if $0<\widetilde{z}$ and $\widetilde{x}<\widetilde{y}$, 
then $(\widetilde{x}\cdot \widetilde{z})<(\widetilde{y}\cdot \widetilde{z})$; 

${\bf \Omega 3}^{\prime }$: if $\widetilde{x}<\widetilde{y}$ then there is 
not, that $\widetilde{y}<\widetilde{x}$ or $\widetilde{x}=\widetilde{y}$ and 
vice versa; 

${\bf \Omega 3}^{\prime \prime }$: for all standard Q-numbers $x$, $y$, $z$: 
$x<y$ or $y<x$ or $x=y$. 

{\bf Proof of the Theorem A.10:} is obvious from Definition A.10 by the 
Theorem A.2. 

{\bf Theorem A.11: }${\bf \Phi }${\bf 8:} If $0<|\widetilde{x}|$ then the 
Q-number $\widetilde{y}$ exists, for which $(\widetilde{x}\cdot \widetilde{y}% 
)=1$.

{\bf Proof of the Theorem A.11:} If $\left\langle x_n\right\rangle \in 
\widetilde{x}$ then by Definition A.10: if 

\[ 
{\bf A}=\left\{ n\in {\bf N}|0<\left| x_n\right| \right\} 
\] 

then ${\bf A}\in {\bf \Phi ix}$. 

In this case: if for the sequence $\left\langle y_n\right\rangle $ : if $% 
n\in {\bf A}$ then $y_n=1/x_n$ 

- then 

\[ 
\left\{ n\in {\bf N}|x_n\cdot y_n=1\right\} \in {\bf \Phi ix}\mbox{.} 
\] 
 
Thus, Q-numbers are fulfilled to all properties of real numbers, except $% 
\Omega $3 \cite{DVS5}. The property $\Omega $3 is accomplished by some weak 
meaning ($\Omega $3' and $\Omega $3''). 

{\bf Definition A.11:} The Q-number $\widetilde{x}$ is {\it the 
infinitesimal Q-number} if the sequence of real numbers $\left\langle 
x_n\right\rangle $ exists, for which: $\left\langle x_n\right\rangle \in 
\widetilde{x}$ and for all positive real numbers $\varepsilon $: 

\[ 
\left\{ n\in {\bf N}||x_n|<\varepsilon \right\} \in {\bf \Phi ix}\mbox{.} 
\] 

Let the set of all infinitesimal Q-numbers be denoted as $I$. 

{\bf Definition A.12:} The Q-numbers $\widetilde{x}$ and $\widetilde{y}$ are 
t{\it he infinite closed Q-numbers} (denote: $\widetilde{x}\approx 
\widetilde{y}$) if $|\widetilde{x}-\widetilde{y}|=0$ or $|\widetilde{x}-% 
\widetilde{y}|$ is infinitesimal. 

{\bf Definition A.13}: The Q-number $\widetilde{x}$ is {\it the infinite 
Q-number} if the sequence $\left\langle r_n\right\rangle $ of real numbers 
exists, for which $\left\langle r_n\right\rangle \in \widetilde{x}$ and for 
every natural number $m$: 

\[ 
\left\{ n\in {\bf N}|m<r_n\right\} \in {\bf \Phi ix}\mbox{.} 
\] 

\subsection{Model} 

Let us define the propositional calculus like to (\cite{MEN63}), but the 
propositional forms shall be marked by the script greek letters. 

{\bf Definition C1: }A set $\Re $ of the propositional forms is{\it \ a 
U-world} if: 

1) if $\alpha _1,\alpha _2,\ldots ,\alpha _n\in \Re $ and $\alpha _1,\alpha 
_2,\ldots ,\alpha _n\vdash \beta $ then $\beta \in \Re $, 

2) for all propositional forms $\alpha $: it is not that $(\alpha \& \left( 
\neg \alpha \right) )\in \Re $, 

3) for every propositional form $\alpha $: $\alpha \in \Re $ or $(\neg 
\alpha )\in \Re $. 

{\bf Definition C2: }The sequences of the propositional forms $\left\langle 
\alpha _n\right\rangle $ and $\left\langle \beta _n\right\rangle $ are {\it % 
Q-equivalent} (denote: $\left\langle \alpha _n\right\rangle \sim 
\left\langle \beta _n\right\rangle $) if 

\[ 
\left\{ n\in {\bf N}|\alpha _n\equiv \beta _n\right\} \in {\bf \Phi ix}% 
\mbox{.} 
\] 

Let us define the notions of {\it the Q-extension of the functions} for like as 
in the Definitions A.5, A.2, A.9, A.5, A.6. 

{\bf Definition C3:} The Q-form $\widetilde{\alpha }$ is {\it Q-real} in the 
U-world $\Re $ if the sequence $\left\langle \alpha _n\right\rangle $ of the 
propositional forms exists, for which: $\left\langle \alpha _n\right\rangle 
\in \widetilde{\alpha }$ and 

\[ 
\left\{ n\in {\bf N}|\alpha _n\in \Re \right\} \in {\bf \Phi ix}\mbox{.} 
\] 

{\bf Definition C4: }The set $\widetilde{\Re }$ of the Q-forms is t{\it he 
Q-extension of the U-world }$\Re $ if $\widetilde{\Re }$ is the set of 
Q-forms $\widetilde{\alpha }$, which are Q-real in $\Re $. 

{\bf Definition C5:} The sequence $\left\langle \widetilde{\Re }% 
_k\right\rangle $ of the Q-extensions is {\it the S-world}. 

{\bf Definition C6: }The Q-form $\widetilde{\alpha }$ is {\it S-real in the 
S-world }$\left\langle \widetilde{\Re }_k\right\rangle $ if 

\[ 
\left\{ k\in {\bf N}|\widetilde{\alpha }\in \widetilde{\Re }_k\right\} \in 
{\bf \Phi ix}\mbox{.} 
\] 

{\bf Definition C7:} The set ${\bf A}$ (${\bf A}\subseteq {\bf N}$) is {\it % 
the regular set} if for every real positive number $\varepsilon $ the 
natural number $n_0$ exists, for which: for all natural numbers $n$ and $m$, 
which are more or equal to $n_0$: 

\[ 
|w_n({\bf A})-w_m({\bf A})|<\varepsilon \mbox{.} 
\] 

{\bf Theorem C1:} If ${\bf A}$ is the regular set and for all real positive $% 
\varepsilon $: 

\[ 
\left\{ k\in {\bf N}|w_k({\bf A})<\varepsilon \right\} \in {\bf \Phi ix}% 
\mbox{.} 
\] 

then 

\[ 
\lim_{k\rightarrow \infty }w_k({\bf A})=0\mbox{.} 
\] 

{\bf Proof of theTheorem C1:} Let be 

\[ 
\lim_{k\rightarrow \infty }w_k({\bf A})\neq 0\mbox{.} 
\] 

That is the real number $\varepsilon _0$ exists, for which: for every 
natural number $n^{\prime }$ the natural number $n$ exists, for which: 

\[ 
n>n^{\prime }\mbox{ and }w_n({\bf A})>\varepsilon _0. 
\] 

Let $\delta _0$ be some positive real number, for which: $\varepsilon 
_0-\delta _0>0$. Because ${\bf A}$ is the regular set then for $\delta _0$ 
the natural number $n_0$ exists, for which: for all natural numbers $n$ and $% 
m$, which are more or equal to $n_0$: 

\[ 
|w_m({\bf A})-w_n({\bf A})|<\delta _0\mbox{.} 
\] 

That is 

\[ 
w_m({\bf A})>w_n({\bf A})-\delta _0\mbox{.} 
\] 

Since $w_n({\bf A})\geq \varepsilon _0$ then $w_m({\bf A})\geq \varepsilon 
_0-\delta _0$. 

Hence, the natural number $n_0$ exists, for which: for all natural numbers $% 
m $: if $m\geq n_0$ then $w_m({\bf A})\geq \varepsilon _0-\delta _0$. 

Therefore, 

\[ 
\left\{ m\in {\bf N}|w_m({\bf A})\geq \varepsilon _0-\delta _0\right\} \in 
{\bf \Phi ix}\mbox{.} 
\] 

and by this Theorem condition: 

\[ 
\left\{ k\in {\bf N}|w_k({\bf A})<\varepsilon _0-\delta _0\right\} \in {\bf % 
\Phi ix}\mbox{.} 
\] 

Hence, 

\[ 
\left\{ k\in {\bf N}|\varepsilon _0-\delta _0<\varepsilon _0-\delta 
_0\right\} \in {\bf \Phi ix}\mbox{.} 
\] 

That is $\emptyset \notin {\bf \Phi ix}$. It is the contradiction for the 
Theorem 2.2. 

{\bf Definition C8:} Let $\left\langle \widetilde{\Re }_k\right\rangle $ be 
a S-world. 

In this case the function ${\mathfrak W}(\widetilde{\beta })$, which has got the 
domain in the set of the Q-forms, has got the range of values in $Q{\bf R}$, 
is defined as the following: 

If ${\mathfrak W}(\widetilde{\beta })=\widetilde{p}$ then the sequence $% 
\left\langle p_n\right\rangle $ of the real numbers exists, for which: $% 
\left\langle p_n\right\rangle \in \widetilde{p}$ and 

\[ 
p_n=w_n\left( \left\{ k\in {\bf N}|\widetilde{\beta }\in \widetilde{\Re }% 
_k\right\} \right) \mbox{.} 
\] 

{\bf Theorem C2:} If $\left\{ k\in {\bf N}|\widetilde{\beta }\in \widetilde{% 
\Re }_k\right\} $ is the regular set and ${\mathfrak W}(\widetilde{\beta }% 
)\approx 1$ then $\widetilde{\beta }$ is S-resl in $\left\langle \widetilde{% 
\Re }_k\right\rangle $. 

{\bf Proof of the Theorem C2: }Since ${\mathfrak W}(\widetilde{\beta })\approx 1$ 
then by Definitions.2.12 and 2.11: for all positive real $\varepsilon $: 

\[ 
\left\{ n\in {\bf N}|w_n\left( \left\{ k\in {\bf N}|\widetilde{\beta }\in 
\widetilde{\Re }_k\right\} \right) >1-\varepsilon \right\} \in {\bf \Phi ix}% 
\mbox{.} 
\] 

Hence, by the point 3 of the Theorem 2.1: for all positive real $\varepsilon 
$: 

\[ 
\left\{ n\in {\bf N}|\left( {\bf N}-w_n\left( \left\{ k\in {\bf N}|% 
\widetilde{\beta }\in \widetilde{\Re }_k\right\} \right) \right) 
<\varepsilon \right\} \in {\bf \Phi ix}\mbox{.} 
\] 

Therefore, by the Theorem C1: 

\[ 
\lim_{n\rightarrow \infty }\left( {\bf N}-w_n\left( \left\{ k\in {\bf N}|% 
\widetilde{\beta }\in \widetilde{\Re }_k\right\} \right) \right) =0\mbox{.} 
\] 

That is: 

\[ 
\lim_{n\rightarrow \infty }w_n\left( \left\{ k\in {\bf N}|\widetilde{\beta }% 
\in \widetilde{\Re }_k\right\} \right) =1\mbox{.} 
\] 

Hence, by Definition.2.3: 

\[ 
\left\{ k\in {\bf N}|\widetilde{\beta }\in \widetilde{\Re }_k\right\} \in 
{\bf \Phi ix}\mbox{.} 
\] 

And by Definition C6: $\widetilde{\beta }$ is S-real in $\left\langle 
\widetilde{\Re }_k\right\rangle $. 

{\bf Theorem C3: }The P-function exists. 

{\bf Proof of the Theorem C3:} By the Theorems C2 and 2.1: ${\mathfrak W}(% 
\widetilde{\beta })$ is the P-function in $\left\langle \widetilde{\Re }% 
_k\right\rangle $.

\end{document}